\documentclass[12pt]{article}   
\usepackage{graphicx,psfrag}
\newtheorem{theorem}{Theorem} [section]
\newtheorem{corollary}[theorem]{Corollary}

\newtheorem{lemma}[theorem]{Lemma}
\newtheorem{prop}[theorem]{Proposition}

\newenvironment {proof} {{\it
Proof.}}{\hspace*{\fill}$\Box$\par\vspace{4mm}}
\usepackage{amsmath}
\usepackage{amsfonts}
\newfont{\bb}{msbm10}

\def\:{\! :\!}

\newcommand\0{{\bf 0}}

\begin{document}

\title{On the $t$-Term Rank of a Matrix}

 \author{Richard A. Brualdi, Kathleen P. Kiernan, \\Seth A. Meyer, Michael W. Schroeder\\
 Department of Mathematics\\
 University of Wisconsin\\
 Madison, WI 53706\\
 {\tt \{brualdi,kiernan,smeyer,schroede\}@math.wisc.edu}
 }

\maketitle

\noindent
{\it Dedicated to Jose Dias da Silva in admiration of his many contributions to linear and multilinear algebra.}

 \begin{abstract} For $t$ a positive integer, the $t$-term rank of a $(0,1)$-matrix $A$ is defined to be the largest number of $1$s in $A$ with at most one 1 in each column and at most $t$ $1$s in each row.
 Thus the $1$-term rank is the ordinary term rank. We generalize some basic results for the term rank to the $t$-term rank, including a formula for the maximum term rank over a nonempty class of $(0,1)$-matrices with the the same row sum and column sum vectors. We also show the surprising result that in such a class there exists a matrix which realizes all of the maximum terms ranks between 1 and $t$.

\medskip
\noindent {\bf Key words and phrases: $t$-term rank, interchange, matrix class, row sum vector, column sum vector } 

\noindent {\bf Mathematics  Subject Classifications: 05A15, 05C50, 05D15.} 
\end{abstract}

\section{Introduction}

Let $A=[a_{ij}]$ be an $m\times n$ matrix which, without loss of generality, we take to be a (0,1)-matrix. The {\it term rank} of    $A$, denoted $\rho(A)$, equals the maximum number of $1$s in $A$ with no two of the $1$s from the same line (row or column). By the well-known K\"onig-Egerv\'ary theorem (see e.g. \cite{BR}, p.~6), $\rho(A)$ equals the minimum number of lines that cover all the
$1$s of $A$:
\[\rho(A)=\min\{e+f: \mbox{$\exists$ a cover of $A$ with $e$ rows and $f$ columns}\}.\]

 If $G$ is the bipartite graph of which $A$ is the bi-adjacency matrix, $\rho(A)$ is the maximum size of a {\it matching} of $G$.

The following generalization of the term rank is motivated by the recent study of combinatorial batch codes (see
\cite{BKMS,PSW} and also \cite{IKOS}). Let $t$ be a positive integer. Then the {\it $t$-term rank} of $A$, denoted $\rho_t(A)$, 
equals the maximum
number of $1$s in $A$ with at most one $1$ in each column and at most $t$ $1$s in each row.  We have $\rho_1(A)=\rho(A)$, and we use both of the notations. For two real $m\times n$ matrices $X=[x_{ij}]$ and $Y=[y_{ij}]$, define $X\le Y$ provided that $x_{ij}\le y_{ij}$ for all $i$ and $j$. Also define \[\sigma(X)=\sum_{i,j}x_{ij}\] the sum of all of the entries of $X$, and 
\[r_i(X)=\sum_{j}x_{ij}\quad  (1\le i\le m),\  s_j(X)=\sum_{i}x_{ij}\quad  (1\le j\le n),\]
the {\it row sums} and {\it column sums} of $X$, respectively.
Then $(r_1(X),\ldots,r_m(X))$ and $(s_1(X),\ldots,s_n(X))$ are, respectively, the {\it row sum vector} and {\it column sum vector} of $X$. It follows that 
\[\rho_t(A)=\max\{\sigma(B):B\le A, r_i(B)\le t\ (1\le i\le m), s_j(B)\le 1\ 
(1\le j\le n)\}.\]

Let $R=(r_1,r_2,\ldots,r_m)$ and $S=(s_1,s_2,\ldots,s_n)$ be two nonnegative integral vectors with $r_1+r_2+\cdots+r_m=s_1+s_2+\ldots+s_n$ . Without loss in generality, we assume  that $R$ and $S$ are monotone non-increasing:
\[r_1\ge r_2\ge\cdots\ge r_m\mbox{ and } s_1\ge s_2\ge\cdots\ge s_n.\]
The class  of all $(0,1)$-matrices with row sum vector $R$ and column sum vector $S$ is denoted by ${\mathcal A}(R,S)$. Note that if $R$ or $S$ has a negative entry, then ${\mathcal A}(R,S)$ is empty. The class ${\mathcal A}(R,S)$ has been heavily  investigated (see \cite{B} for a detailed treatment). In particular, the Gale-Ryser theorem gives necessary and sufficient conditions for 
such a class ${\mathcal A}(R,S)$ to be nonempty. We do not need to use these conditions so we do not state them  here. An alternative criterion for the nonemptiness of ${\mathcal A}(R,S)$ is due to Ford and Fulkerson, and we now review this for later use. Let an $(m+1)\times (n+1)$ matrix
$T=[t_{ij}]$ (in these matrices, rows are indexed by $0,1,\ldots,m$ and columns are indexed by $0,1,
\ldots,n$) be defined by
\begin{equation}\label{eq:FF}
t_{kl}=kl-\sum_{j=1}^l s_j+\sum_{i=k+1}^mr_i\ (0\le k\le m, 0\le j\le n).\end{equation}
It is straightforward to 
check that if a matrix $A\in {\mathcal A}(R,S)$ is partitioned as
\[A=\left[\begin{array}{cc}
X&A_{12}\\ A_{21}&Y\end{array}\right]\mbox{ where $X$ is $k\times l$},\] then
\[t_{kl}=(kl-\sigma(X))+\sigma(Y), \]
the number of $0$s in $X$ plus the number of $1$s in $Y$.
We then have that {\it ${\mathcal A}(R,S)\ne \emptyset$ if and only if
each entry of $T$ is  nonnegative}. In fact, the nonnegativity of $T$ is sufficient for the nonemptiness of ${\mathcal A}(R,S)$ under less restrictive monotonicity assumptions on $R$ and $S$. We review this here as we shall make important use of it later.
An integral vector $U=(u_1,u_2,\ldots,u_n)$ is {\it nearly nonincreasing} provided
\[u_i\ge u_j-1\ (1\le i<j\le n).\]
For example, $(3,4,3,4)$ is nearly nonincreasing. It follows that an integral vector is nearly nonincreasing if and only if it can be made into a monotone nonincreasing vector by the addition of some $(0,1)$-vector. The following theorem is stated and proved as Theorem 2.1.4 in \cite{B}.

\begin{theorem} \label{th:FF}
Let $R=(r_1,r_2,\ldots,r_m)$ and $S=(s_1,s_2,\ldots,s_n)$ be nonnegative integral vectors such that $S$ is nearly nonincreasing and $r_1+r_2+\cdots+r_m=s_1+s_2+\cdots+s_n$. Then ${\mathcal A}(R,S)$ is nonempty if and only if
\[t_{kl}\ge 0\quad (0\le k\le m, 0\le l\le n),\]
where $t_{kl}$ is defined as in $(\ref{eq:FF})$.
\end{theorem}

In Section 2, we determine some basic properties of the $t$-term rank, in particular, some basic properties of the $t$-term rank over the matrices in a class ${\cal A}(R,S)$.
In Section 3, we show that the $t$-term rank of a semi-regular matrix depends only on $t$ and the dimensions of the matrix.
In Section 4, we obtain a formula for the maximum $t$-term rank over a class ${\mathcal A}(R,S)$. 
In Section 5, we prove a surprising theorem concerning the maximum $t$-term rank in a nonempty class ${\mathcal A}(R,S)$; we show, in particular, that there is a matrix in ${\mathcal A}(R,S)$ which attains the maximum 
$t'$-term rank for all $t'$ with $1\le t'\le t$.

\section{The $t$-term rank}

In this section we discuss some elementary properties of $\rho_t(A)$.

Let $A$ be an $m\times n$ $(0,1)$-matrix and let $t$ be a positive integer.
Let $A^{(t)}$ be the $tm\times n$ matrix obtained by stacking up $t$ copies of $A$. Thus, for instance, 
\[A^{(3)}=\left[\begin{array}{c}
A\\ \hline
A\\ \hline
A\end{array}\right].\]
We immediately get
\[\rho_t(A)=\rho(A^{(t)}).\]

\begin{prop}\label{prop:KE}
Let $A$ be an $m\times n$ $(0,1)$-matrix and let $t$ be a positive integer. Then
\[\rho_t(A)=\min\{te+f: \mbox{$\exists$ a cover of $A$ with $e$ rows and $f$ columns}\}.\]
\end{prop}

\begin{proof} This follows by applying the K\"onig-Egerv\'ary theorem to $A^{(t)}$ noting that a column used in a cover of $A^{(t)}$ corresponds to the same column in each of the copies of $A$ making up $A^{(t)}$.\end{proof}

\noindent
{\bf Remarks:} (1) Let $A$ be an $m\times n$ $(0,1)$-matrix and, in the standard way, regard $A$ as the incidence matrix of a family ${\mathcal C}=(C_1,C_2,\ldots,C_m)$ of subsets of a set $X=\{x_1,x_2,\ldots,x_n\}$ of $n$ elements. Then P.~Hall's theorem  for the existence of a system of distinct representatives of $\mathcal C$ can be used to give a different but equivalent expression for the term rank of $A$, namely,
\[\rho(A)=\min_{K\subseteq \{1,2,\ldots,m\}}
\{|\cup_{i\in K}C_i| +(m-|K|)\}.\]
In general, we obtain  a different but equivalent expression for $\rho_t(A)$, namely,
\[\rho_t(A)=\min_{K\subseteq\{1,2,\ldots,m\}}
\{|\cup_{i\in K}C_i| +t(m-|K|)\}.\]

\smallskip\noindent
(2) Let $A$ be an $m\times n$  $(0,1)$-matrix. For $K\subseteq \{1,2,\ldots,n\}$, let $A[*,K]$ denote the $m\times |K|$ submatrix of $A$ determined by the columns with index in $K$. Then the set of all such $K$ with $\rho(A[*,K])=|K|$ are the independent sets of a matroid 
${\bf M}(A)$ on $\{1,2,\ldots,n\}$. Such matroids are {\it transversal matroids}, and the rank of ${\bf M}(A)$ is  $\rho(A)$ (see e.g. \cite{RAB,Ox}). Let $t$ be a positive integer.  The transversal matroid ${\bf M}(A^{(t)})$ has rank $\rho_t(A)$  and is the matroid union (see e.g. \cite{Ox}) of $t$ copies of ${\bf M}(A)$:
\[{\bf M}(A^{(t)})={\bf M}(A)\vee {\bf M}(A)\vee\cdots\vee {\bf M}(A)\quad
(t \mbox{ copies of } {\bf M}(A)).\]
The {\it bases} (maximal independent sets) of ${\bf M}(A^{(t)})$ are those sets $K\subseteq \{1,2,\ldots,n\}$ with $|K|=\rho_t(A)$ which can be partitioned into sets $K_1,K_2,\ldots,K_t$ such that 
$A[*,K_i]$ has at most one $1$ in each of its rows and exactly one 1 in each of its columns ($i=1,2,\ldots,t)$.
Because we are dealing with matroids, there is always such a basis 
of ${\bf M}(A^{(t)})$  such that $K_1\cup K_2\cup\cdots\cup K_j$ is a basis
of ${\bf M}(A^{(j)})$ for each $j=1,2,\ldots,t-1$. This remark establishes the following proposition.

\begin{prop}\label{prop:new}
Let $A$ be an $m\times n$  $(0,1)$-matrix, and let $t\ge 2$ be an  integer. Then there exists a $(0,1)$-matrix $B\le A$  where
\begin{enumerate}
\item[\rm (i)] 
$r_i(B)\le t$\ $(i=1,2,\ldots,m)$ and $s_j(B)\le 1$\  $(j=1,2,\ldots,n)$.
\item[\rm (ii)] $\rho_t(A)=\sigma(B)=\sigma(B[*,K])=|K|$ where $K=\{j: 1\le j\le n, s_j=1\}$.
\item[\rm (iii)] there is a $(0,1)$-matrix $C\le B$ such that
$r_i(C)\le t-1$\ $(i=1,2,\ldots,m)$,  $s_j(C)\le 1$\  $(j=1,2,\ldots,n)$,
and $\rho_{t-1}(A)=\sigma(C)=\sigma(C[*,K'])$ for some $K'\subseteq K$ with $|K'|=\rho_{t-1}(A)$.
\end{enumerate}
\end{prop}

Let $A$ be  an $m\times n$ $(0,1)$-matrix with  $m\le n$ and with at least one 1 in each  column.
Then clearly $\rho_t(A)$ as a function of $t$ strictly increases until it takes the value $n$.
We define the {\it strength} of the $m\times n$ $(0,1)$-matrix $A$, denoted by $\gamma(A)$, to be the smallest positive integer $t$ such that $\rho_t(A)=n$, that is,  the smallest positive integer $t$ such that there exists an $m\times n$ $(0,1)$-matrix $B\le A$ which has exactly one  1 in every column and at most $t$ $1$s in every row. Since $A$ has at least one 1 in each column, $\gamma(A)$ is well-defined.
 If we define $\rho_0(A)=0$, we have
\[\rho_0(A)<\rho_1(A)<\rho_2(A)<\cdots< \rho_{\gamma(A)-1}(A)<\rho_{\gamma(A)}=n.\]
 The strength $\gamma(A)$ equals the smallest integer $t$ such that $\{1,2,\ldots,n\}$ is an independent set of the matroid ${\bf M}(A^{(t)})$. 
It follows that
\[\gamma(A)\le \max\{r_1(A),r_2(A),\ldots,r_m(A)\}.\]

\begin{prop}
\label{prop:concave}
If $A$ is an $m\times n$ $(0,1)$-matrix, the sequence
$\rho_0(A),\rho_1(A),\rho_2(A),\ldots$ satisfies
\[\rho_k(A)-\rho_{k-1}(A)\ge \rho_{k+1}(A)-\rho_k(A)\ (k\ge 1).\]
\end{prop}

\begin{proof}
This proposition  is an easy  consequence of the matroidal connections already discussed. 
There is a basis $K=K_1\cup K_2\cup \cdots\cup K_{k+1}$
of ${\bf M}(A^{(k+1)})$  such that $K_1\cup K_2\cup\cdots\cup K_{j}$ is a basis
of ${\bf M}(A^{(j)})$ for each $j=1,2,\ldots,k$.
This implies that
$\rho_{k+1}(A)=\rho_{k}(A)+|K_{k+1}|$. Similarly,
$\rho_{k}(A)=\rho_{k-1}(A)+|K_k|$. Since clearly $|K_k|\ge |K_{k+1}|$, the proposition follows.
\end{proof}

A basic property of a nonempty class ${\mathcal A}(R,S)$ is that starting from any one matrix $A\in {\mathcal A}(R,S)$, we can get to any other matrix by a sequence of interchanges where an {\it interchange}  replaces
one of the $2\times 2$ submatrices 
\[\left[\begin{array}{cc} 1&0\\0&1\end{array}\right]\mbox{ and }
\left[\begin{array}{cc} 0&1\\1&0\end{array}\right]\]
with the other (see e.g. \cite{B}). A single interchange can change the term rank by 
at most 1, either positively or negatively \cite{B}.
We show that a similar conclusion holds for the $t$-term rank in general.

\begin{prop}\label{prop:interchange}
Let $A$ be an $m\times n$ $(0,1)$-matrix, and let $t$ be a positive integer. Let $A'$ be obtained from $A$ by a single interchange.
Then
\[\rho_t(A)-1\le \rho_t(A')\le \rho_t(A)+1.\]
\end{prop}

\begin{proof}
Consider the matrix $A^{(t)}$. Then there is a cover of $A^{(t)}$ with
$\rho_t(A)$ lines with the property that if row $i$ is used in one copy of $A$, then it is used in every copy of $A$ (if not then we could eliminate row $i$ in each copy of $A$ in which it is used). Thus in this cover, the same rows and columns are used in each copy of $A$. Outside the union of these rows and columns there is a zero matrix in each copy of $A$. The matrix $A'$ can have at most one 1 in its positions corresponding to this zero matrix. Hence the matrix $A'^{(t)}$
can be covered by using one additional column and no additional rows, and so
\[\rho_t(A')=\rho(A'^{(t)})\le \rho(A^{(t)})+1=\rho_t(A)+1.\]
Since an interchange is reversible, the inequality
$\rho_t(A)-1\le \rho_t(A)$ also follows.
\end{proof}

\smallskip\noindent
{\bf Example:} Consider the $7\times 9$ matrix
\[A=\left[\begin{array}{ccccccccc}
0&0&0&1&0&0&1&0&0\\
0&0&0&0&1&0&1&0&0\\
0&0&0&1&0&1&0&1&1\\
1&0&1&0&0&0&0&0&0\\
0&1&0&0&0&0&0&0&0\\
0&0&1&0&0&0&0&0&0\\
1&0&0&0&0&0&0&0&0\end{array}\right],\]
where it is straightforward to check that $\rho_1(A)=6$ and $\rho_2(A)=8$.
The interchange using the $2\times 2$ submatrix in
rows and columns 3 and 4 produces the matrix
\[A'=\left[\begin{array}{ccccccccc}
0&0&0&1&0&0&1&0&0\\
0&0&0&0&1&0&1&0&0\\
0&0&1&0&0&1&0&1&1\\
1&0&0&1&0&0&0&0&0\\
0&1&0&0&0&0&0&0&0\\
0&0&1&0&0&0&0&0&0\\
1&0&0&0&0&0&0&0&0\end{array}\right],\]
where $\rho_1(A')=7$ but $\rho_2(A')=\rho_2(A)=8$. Thus even though $\rho_1(A)$ increases, $\rho_2(A)$ does not, although there is room for an increase since $\rho_2(A)=8<9$.
\smallskip

We now show that a single interchange which increases  the $(t-1)$-term rank cannot decrease the $t$-term rank (cf. Proposition \ref{prop:interchange}). As shown in the preceding example, the $t$-term rank may not change.

\begin{prop}\label{prop:single}
Let  $A'$ be obtained  from a matrix $A$ by a single interchange, and let $t\ge 2$ be an   integer.
Assume that $\rho_{t-1}(A')=\rho_{t-1}(A)+1$. Then
\[\rho_{t}(A)\le \rho_{t}(A')\le \rho_{t}(A)+1.\]
\end{prop}

\begin{proof}
By Proposition \ref{prop:interchange} we need only show that $\rho_t(A')\ge \rho_t(A)$.  Let $B\le A$ and $C\le A$ be  matrices whose existences are established in Proposition \ref{prop:new} whose notation we now use. The interchange that produces $A'$ from $A$ and increases $\rho_{t-1}(A)$
must take place in one column of $K'$ and one column from the complement $\overline{K'}$ of $K'$. Hence the matrix $C$ can lose at most one 1 as a result of the interchange. But the new 1 now in a column of $\overline{K'}$ either replaces a 1 of $B$ (when the other column of the interchange is in $K
\setminus K'$) or can be used to add a new 1 to $B$ (when the other column of the interchange is in $\overline{K}$). Hence $\rho_t(A)$ does not decrease.
\end{proof}

\smallskip\noindent
{\bf Example:} 
Let
\[A=\left[\begin{array}{cccccc}
1&0&0&0&0&1\\
0&1&0&0&0&1\\
0&0&1&1&1&0\end{array}\right],\]
Then $\rho_1(A)=3$ and $\rho_2(A)=5$. Applying the interchange in the lower right corner of $A$ gives the matrix
\[A'=\left[\begin{array}{cccccc}
1&0&0&0&0&1\\
0&1&0&0&1&0\\
0&0&1&1&0&1\end{array}\right].\]
where $\rho_1(A')=3$ and $\rho_2(A')=6$.  Both $\rho_1$ and $\rho_2$ can increase by 1 after an interchange, as the example
\[\left[\begin{array}{cccc}
0&1&1&1\\
1&0&0&0\\
1&0&0&0\\
1&0&0&0\end{array}\right]\rightarrow
\left[\begin{array}{cccc}
1&0&1&1\\
0&1&0&0\\
1&0&0&0\\
1&0&0&0\end{array}\right]\]
shows, where $\rho_1(A)=2,\rho_2(A)=3, \rho_1(A')=3,\mbox{ and }
\rho_2(A')=4$.

\section{Semiregular Classes}

Let $m$ and $n$ be positive integers, and let $k$ and 
$l$ be positive integers such that $km=nl$.  Then ${\mathcal A}(m,n;k,l)$
denotes the class of all $m\times n$ $(0,1)$-matrices with $k$\ $1$s in each row and $l$\ $1$s in each column. Since $km=ln$, this class is nonempty.  Matrices in a class ${\mathcal A}(m,n;k,l)$ are called
{\it semiregular}; in case $m=n$, and thus $k=l$, the matrices are {\it regular}. We show that for each positive integer $t$, the $t$-term rank is constant on ${\mathcal A}(m,n;k,l)$.

\begin{theorem}\label{th:semi} For a nonempty class ${\mathcal A}(m,n;k,l)$ and $t$ a positive integer, we have
\[\rho_t(A)=\min\{tm,n\} \mbox{ for all }A\in {\mathcal A}(m,n;k,l).\]
\end{theorem}

\begin{proof} Let $A\in {\mathcal A}(m,n;k,l)$.
First consider the case where $tm\le n$ so that $t\le n/m$.  We need 
to show that $\rho_t(A)=tm$.  Suppose  that $\rho_t(A)<tm$ so that by Proposition \ref{prop:KE} there exist $e$ rows and $f$ columns of $A$ which cover all the $1$s of $A$ where $te+f<tm$.
Thus after row and column permutations, we can take $A$ in the form
\begin{equation}\label{eq:form}
\left[\begin{array}{cc}
A_1&A_{12}\\
A_{21}&O\end{array}\right] \ (A_1 \mbox{ is } e\times f).\end{equation}
 We then have 
 \[f<t(m-e)\le \frac{n}{m}(m-e)=\frac{k}{l}(m-e),\]
 and so
 \[fl<k(m-e).\]
 This implies that the total number of $1$s in $A_{21}$ is strictly greater than the total number of $1$s in the first $f$ columns of $A$, a contradiction.

Now consider the case where $n\le tm$ so that $t\ge n/m$. Suppose that $\rho_t(A)<n$
so that we may assume that $A$ has the form (\ref{eq:form}) where
$te+f<n$. We then have
\[\frac{n}{m}e+f\le te+f<n\mbox{ 
and hence } ne<m(n-f).\]
Since $n=(km)/l$, this
gives
\[ke<l(n-f),\]
implying that the total number of $1$s in $A_{12}$ is strictly greater than the total number of $1$s in the first $e$ rows of $A$, a contradiction. 
\end{proof}

\begin{corollary}\label{cor:sem}
For a nonempty class ${\mathcal A}(m,n;k,l)$ with $m\le n$, we have
\[\gamma(A)=\left\lceil\frac{n}{m}\right\rceil \mbox{ for all } A\in {\mathcal A}(m,n;k,l).\]
\end{corollary}

\section{Formula for Maximum $t$-Term Rank}

A formula of Ryser (see \cite{B}, p.~71) gives the following formula for the maximum term rank, denoted  by $\overline{\rho}(R,S)$ or $\overline{\rho_1} (R,S)$, of matrices in a nonempty class ${\mathcal A}(R,S)$. Assuming that $R$ and $S$ are monotone nonincreasing,  we have
\begin{equation}\label{eq:ryser}
\overline{\rho}(R,S)=\min\{t_{ef}+e+f: 0\le e\le m, 0\le f\le n\}.\end{equation}
Our first goal in this section is to generalize this formula to the {\it maximum $t$-term rank}, denoted  by $\overline{\rho_t}(R,S)$, of matrices in ${\mathcal A}(R,S)$.
To do this, we generalize the proof given by Brualdi and Ross (see \cite{B}, p.~69--71) for $\overline{\rho}(R,S)$. We shall make use of the following existence theorem which in  the general form given is due to Anstee
(see \cite{B}, p.~189). We state it in the transposed form for our purposes.

\begin{theorem} Let $R=(r_1,r_2,\ldots,r_m)$ and $S=(s_1,s_2,\ldots,s_n)$ be nonnegative integral vectors.
Let $k$ be a nonnegative integer, and let $(k_1,k_2,\ldots,k_n)$ be a prescribed vector of integers with $k\le k_i\le k+1$ for $i=1,2,\ldots,n$. 
Let $R'=(r_1',r_2',\ldots,r_m')$ where $r_i'\le r_i$ for $i=1,2,\ldots,m$, and let $S'=(s_1',s_2',\ldots,s_n')=(s_1-k_1,s_2-k_2,\ldots,s_n-k_n)$. 
Then there exist matrices $A\in {\mathcal A}(R,S)$ 
and $A'\in {\mathcal A}(R',S')$ with $A'\le A$ if and only if both of the classes ${\mathcal A}(R,S)$ and ${\mathcal A}(R',S')$ are nonempty.
\end{theorem}

We shall apply this theorem in the following form.

\begin{corollary}\label{cor:class}
Let $R=(r_1,r_2,\ldots,r_m)$ and $S=(s_1,s_2,\ldots,s_n)$ be nonnegative integral vectors. Let $t$ be a positive integer, and let
$R'=(r_1',r_2',\ldots,r_m')$ be an integral vector satisfying $r_i-t\le r_i'\le r_i$ for $i=1,2,\ldots,m$. Let $S'=(s_1',s_2',\ldots,s_n')$ 
be an integral vector satisfying $s_j-1\le s_j'\le s_j$ for $j=1,2,\ldots,n$.
Then there exist matrices $A\in {\mathcal A}(R,S)$ 
and $A'\in {\mathcal A}(R',S')$ with $A'\le A$ if and only if both of the classes ${\mathcal A}(R,S)$ and ${\mathcal A}(R',S')$ are nonempty.
\end{corollary}

Note that for any such $R'$ as given in the corollary, the matrices $A$ and $A'$ in Corollary \ref{cor:class} satisfy
$A-A'\le A$ is a $(0,1)$-matrix with at most $t$ $1$s in each row and at most one 1 in each column. Thus $\rho_t(A)\ge \sum_{j=1}^n (s_j-s_j')$.

The formula generalizing Ryser's formula (\ref{eq:ryser}) is contained in the following theorem. We now write $t_{ef}(R,S)$ in place of $t_{ef}$, since we will have to calculate these numbers for different row sum and column sum vectors.

\begin{theorem}\label{th:maxterm}
Let $R=(r_1,r_2,\ldots,r_m)$ and $S=(s_1,s_2,\ldots,s_n)$ be monotone nonincreasing,  nonnegative integral vectors such that ${\mathcal A}(R,S)$ is nonempty. Let $t$ be a positive integer. Then
\begin{equation}\label{eq:rysergen}
\overline{\rho_t}(R,S)=\min\{t_{ef}(R,S)+te+f: 0\le e\le m, 0\le f\le n\}.\end{equation}
\end{theorem}

\begin{proof}
The proof will be given using several lemmas.
In the first lemma we consider the extreme case where $n=tm$
and $\overline{\rho_t} (R,S)=tm$, the largest it could possibly be. 

\begin{lemma}\label{lem:main1} Let $t$ be a positive integer.
Let $R=(r_1,r_2,\ldots,r_m)$ and $S=(s_1,s_2,\ldots,s_{tm})$ be monotone nonincreasing,  nonnegative integral vectors such that ${\mathcal A}(R,S)$ is nonempty.  Then there exists a matrix $A\in {\mathcal A}(R,S)$ with $\rho_t(A)=tm$ if and only if 
\[t_{ef}(R,S)+te+f \ge tm\quad  (0\le e\le m, 0\le f\le tm).\]
\end{lemma}

\begin{proof} Let $R'=(r_1-t,r_2-t,\ldots,r_m-t)$ and $S'=(s_1-1,s_2-1,\ldots, s_{tm}-1)$. 
It follows from Corollary \ref{cor:class} that the desired matrix $A$ exists if and only if the class ${\mathcal A}(R',S')$ is also nonempty.
Applying Theorem \ref{th:FF} to ${\cal A}(R',S')$,
we see that ${\mathcal A}(R',S')\ne\emptyset$ if and only if $t_{ef}(R',S')\ge 0$ for all $e$ and $f$. An easy calculation shows that
\begin{eqnarray*}
t_{ef}(R',S')&=&t_{ef}(R,S)+f-t(m-e)\\
&=&t_{ef}(R,S)+te+f-tm\quad (0\le e\le m,0\le f\le tm),\end{eqnarray*}
and the lemma follows.
\end{proof}

Thus in the case where $\overline{\rho_t}$ is achieved with $t$ $1$s in every row, we have the following.

\begin{corollary}\label{cor:main1} 
With the hypotheses in Lemma $\ref{lem:main1}$, 
$\overline{\rho_t}(R,S)=tm$  if and only if 
\[\min\{t_{ef}(R,S)+te+f: 0\le e\le m, 0\le f\le tm\}=tm.\]
\end{corollary}

We now show that we can reduce the evaluation of $\overline{\rho_t}(R,S)$ to the situation where Lemma \ref{lem:main1} applies.

\begin{lemma}\label{lem:main2} Let $t$\ be a positive integer.
Let $A$ be an $m\times n$ $(0,1)$-matrix and let $p$ be an integer
with $0\le p\le tm$. Let $l$ be the (smallest) nonnegative integer such that
$tl\ge n-p$, and let $q=tl-(n-p)$.
 Let $A^*$ be the $(m+l)\times t(m+l)$ $(0,1)$-matrix defined by
\[A^*=\left[\begin{array}{ccc}
O_{l,tm-p}&J_{l,n}&J_{l,q}\\
J_{m,tm-p}&A&O_{m,q}
\end{array}\right]\]
where $J$ and $O$ denote matrices of all $1$s and of all $0$s, respectively, of the indicated sizes. 
Then $\rho_t(A)\ge p$ if and only if $\rho_t(A^*)=t(m+l).$
\end{lemma}

\begin{proof}
First suppose that $\rho_t(A)\ge p$. Thus there exists a $(0,1)$-matrix
$B\le A$ with at most $t$ $1$s in each row and at most one $1$ in each column, such that $\sigma(B)=p$.  For any row sum $r$ of $B$,  $tm-p\ge t-r$ if and only if $r\ge p-t(m-1)$, and the latter surely holds since each row sum of $B$ is at most $t$. Thus $tm-p\ge t-r$, and it follows that there exists a $(0,1)$-matrix $C\le J_{m,tm-p}$
such that  the $m$ by $(tm-p+n)$ matrix $[C\; B]$ has exactly $t$ $1$s in each row and at most one 1 in each column with those columns not containing 1s coming from $B$.

There exists a $(0,1)$-matrix $E\le J_{ln}$ with at most $t$ 1s in each row and, since $q\le tl$,  a $(0,1)$-matrix $D\le J_{lq}$ with at most $t$ 1s in each row and exactly one 1 in each column,  such that the matrix
\[T=\left[\begin{array}{ccc}
O_{l,tm-p}&E&D\\
C&B&O_{m,q}
\end{array}\right]\]
satisfies $T\le  A^*$ and has exactly $t$ $1$s in each row and exactly one $1$ in each column; thus $\rho_t(A^*)=t(m+l)$.

Conversely, if $\rho_t(A^*)=t(m+l)$, then there exist $tm$ $1$s from the last $m$ rows of $A^*$ with $t$ $1$s from each row and at most one $1$ from each column. At most $tm-p$ of these $1$s come from
$J_{m,tm-p}$ and hence at least $p$ come from $A$. Hence $\rho_t(A)\ge p$.
\end{proof}

We remark that in the definition of $A^*$, we may replace $O_{l,tm-p}$ with $J_{l,tm-p}$ and Lemma \ref{lem:main2} continues to hold. The reason for doing so now is that the row and column sums of 
\begin{equation}\label{eq:newA*}
A^*=\left[\begin{array}{ccc}
J_{l,tm-p}&J_{l,n}&J_{l,q}\\
J_{m,tm-p}&A&O_{m,q}
\end{array}\right]\end{equation}
are monotone nonincreasing. In applying Lemma \ref{lem:main2}, we shall use $A^*$ as given in (\ref{eq:newA*}).  Let $A$ be a matrix in a class ${\mathcal A}(R,S)$. Then, where $R^*$ and $S^*$ are the row sum and column sum vectors, respectively, of the matrix $A^*$ in (\ref{eq:newA*}),  the matrix $A^*$ defines a class ${\mathcal A}(R^*,S^*)$ with monotone nonincreasing row and column sum vectors $R^*$ and $S^*$ which we make use of below. The matrices in ${\mathcal A}(R^*,S^*)$ are  exactly those matrices (\ref{eq:newA*}) obtained as $A$ varies over the class ${\mathcal A}(R,S)$.

\begin{lemma}\label{lem:main3} Let $t$ be a positive integer, and 
let $R=(r_1,r_2,\ldots,r_m)$ and $S=(s_1,s_2,\ldots,s_n)$ be monotone nonincreasing, nonnegative integral vectors such that the class ${\mathcal A}(R,S)$ is nonempty. Let $p$ be an integer with $0\le p\le tm$. Then there exists a matrix $A\in {\mathcal A}(R,S)$ such that $\rho_t(A)\ge p$ if and only if
\begin{equation}\label{eq:ineq}
t_{ef}(R,S)+te+f\ge p\quad (0\le e\le m, 0\le f\le n).\end{equation}
\end{lemma}

\begin{proof}
First suppose that there is a matrix $A\in {\mathcal A}(R,S)$ such that $\rho_t(A)\ge p$. Let $e$ and $f$ be integers with $0\le e\le m, 0\le f\le n$.  If $te+f\ge p$, then (\ref{eq:ineq}) holds. Now suppose that
$te+f\le p$. Then it follows from Proposition \ref{prop:KE} that the $t$-term rank of the matrix obtained from $A$ by replacing 
its lower right $(m-e)\times (n-f)$ submatrix $Y$ with all 0s is at most
$te+f$. Since $\rho_t(A)\ge p$, $Y$ must have at least $p-(te+f)$ $1$s. Thus $t_{ef}(R,S)\ge p-(te+f)$ and so (\ref{eq:ineq}) holds.

 Now assume that
(\ref{eq:ineq}) holds. We need to show the existence of a matrix $A\in {\mathcal A}(R,S)$ with $\rho_t(A)\ge p$.
By Lemmas \ref{lem:main1} and  \ref{lem:main2}, there exists a matrix $A\in {\mathcal A}(R,S)$ with $\rho_t(A)\ge p$ if and only if 
\begin{equation}\label{eq:needed}
t_{ef}(R^*,S^*)+te+f\ge t(m+l)\quad (0\le e\le m+l,0\le f\le t(m+l)).
\end{equation}

If $e= m+l$ or $f= t(m+l)$, then (\ref{eq:needed}) holds. 
Thus we may assume that $e<m+l$ and $f<t(m+l)$.
Six cases need to be considered, according to which of the six matrices 
in the definition (\ref{eq:newA*}) contains the $(e,f)$-entry of $A^*$.
(Keep in mind that the row and column indices of $T(R^*)=[t_{ij}(R^*,S^*)]$ start with 0.)

\smallskip\noindent
Case 1: $0\le e< l$ and $0\le f\le tm-p$. We have, since $e<l$,
\begin{eqnarray*}
t_{ef}(R^*,S^*)+te+f&\ge&t_{00}(R,S)+ (tm-p-f)+n+q+ te+f\\
&\ge&p+ (tm+n+q-p)+te\\
&=& tm+te+n+q\\
&\ge &tm+tl.\end{eqnarray*}

\smallskip\noindent
Case 2: $0\le e< l$ and $tm-p\le f<tm-p+n$. 
This case is similar to Case 1.

\smallskip
\noindent
Case 3: $0\le e<l$ and $tm-p+n\le f<t(m+l)$. This case is  also similar to Case 1.

\smallskip\noindent
Case 4: $l\le e< m+l$ and $0\le f<tm-p$. Again this case is similar to Case 1.

\smallskip\noindent
Case 5: $l\le e<m+l$ and $tm-p\le  f<tm-p+n$.  Then
\begin{eqnarray*}
t_{ef}(R^*,S^*)+te+f &\ge & t_{e-l,f-(tm-p)}(R,S)+te+f\\
&=& t_{e-l,f-(tm-p)}(R,S)+tl+t(e-l)+(tm-p)+(f-(tm-p))\\
&\ge & p+tl+(tm-p)\\
&=&t(m+l).\\
\end{eqnarray*}

\smallskip\noindent
Case 6: $l\le e< m+l$ and $ tm-p+n\le f<t(m+l)$. A last case similar to Case 1.
\end{proof}

The proof of Theorem \ref{th:maxterm} now follows from Lemma \ref{lem:main3}.

\end{proof}

\section{Joint Realization of $t$-Term Ranks}

A theorem of Haber (see \cite{B}, p.~69) asserts that there exists a matrix $A$ in a nonempty class ${\mathcal A}(R,S)$ with maximum term rank $\overline{\rho}=\overline{\rho}(R,S)$ having $1$s in positions
$(1,\overline{\rho}), (2,\overline{\rho}-1),\ldots,(\overline{\rho},1)$; in particular, the leading $\overline{\rho}\times\overline{\rho}$ submatrix of $A$ has term rank equal to $\overline{\rho}$. In this section, using Theorem \ref{th:FF}, we obtain a significant extension  of this result, the proof of which is based on the following lemma.

\begin{lemma}\label{lem:joint}
Let $R$ and $S$ be monotone nonincreasing nonnegative integral vectors such that ${\mathcal A}(R,S)$ is nonempty. Let $t$ be a positive integer such that $\overline{\rho_t}=\overline{\rho_t}(R,S)\le n$. Consider the integer partition $R'=(r'_1,r'_2,\ldots,r'_m)$ of the integer $\overline{\rho_t}$
defined by
\[R'=(\underbrace{t,t,\ldots,t}_{p_t}, \underbrace{t-1,t-1,\ldots,t-1}_{p_{t-1}},\ldots,\underbrace{1,1,\ldots,1}_{p_1},\underbrace{0,0,\ldots,0}_{p_0})\]
where
\begin{equation}\label{eq:last}
p_k=\left\{\begin{array}{ll}
\overline{\rho_t}-\overline{\rho_{t-1}}&\mbox{ if $k=t$,}\\
2\overline{\rho_k}-\overline{\rho_{k+1}}-\overline{\rho_{k-1}}&\mbox{ if $1\le k<t$,}\\
m-\overline{\rho_1}&\mbox{ if $k=0$.}\end{array}\right.\end{equation}
Let $S'=(s'_1,s'_2,\ldots,s'_n)$ where 
\[s'_j=\left\{\begin{array}{ll} 1&\mbox{ if $1\le j\le \overline{\rho_t}$,}
\\
0&\mbox{ if $\overline{\rho_t}<j\le n$.}\end{array}\right.\]
Then there exists a matrix $A\in {\mathcal A}(R,S)$ and a matrix $C\in {\mathcal A}(R',S')$ such that $C\le A$. 
\end{lemma}

\begin{proof} First we note that by Proposition \ref{prop:concave}, the integers $p_k$ are all nonnegative.
Let $R''=R-R'=(r_1'',r_2'',\ldots,r_m'')$ and $S''=S-S'=(s_1'',s_2'',\ldots,s_n'')$. 
By hypothesis ${\mathcal A}(R,S)$ is nonempty, and so it suffices by Corollary \ref{cor:class} to show that ${\mathcal A}(R'',S'')$ is nonempty.
Since $S'$ is montone nonincreasing, $S''$ is nearly nonincreasing.
Since $r_1''+r_2''+\cdots+r_m''=s_1''+s_2''+\cdots+s_n''$. it now follows by Theorem \ref{th:FF} that we need only verify that
\[t_{ef}(R'',S'')\ge 0\quad  (0\le e\le m, 0\le f\le n).\]
We note that since ${\mathcal A}(R,S)\ne\emptyset$, we have that $t_{ef}(R,S)\ge 0$ for all $e$ and $f$.

We calculate that
\begin{eqnarray*}
t_{ef}(R'',S'')&=&ef+\sum_{i=e+1}^mr''_i-\sum_{j=1}^f s''_j\\
&=&ef+\sum_{i=e+1}^m  r_i-\sum_{i=e+1}^m r'_i-\sum_{j=1}^fs_j+\min\{f,\overline{\rho_t}\}\\
&=& t_{ef}(R,S)-\sum_{i=e+1}^m r'_i+\min\{f,\overline{\rho_t}\}.\\
\end{eqnarray*}

If $f\ge \overline{\rho_t}$, then 
\[
-\sum_{i=e+1}^m r'_i +\min\{f,\overline{\rho_t}\} \ge -\sum_{i=1}^m r'_i +\overline{\rho_t}=0;\]
hence, $t_{ef}(R'',S'')\ge t_{ef}(R,S)\ge 0$. Also, if $e\ge \overline{\rho_1}$, then $-\sum_{i=e+1}^m r'_i =0$ and  again $t_{ef}(R'',S'')\ge 0$.

 Now assume that 
$f<\overline{\rho_t}$  and $e<\overline{\rho_1}$. We now use Theorem \ref{th:maxterm} which implies that for $1\le k\le t$, 
\begin{equation}\label{eq:end}
\overline{\rho_k}\le t_{ef}(R,S)+ke+f\quad (0\le e\le m,0\le f\le n).\end{equation}
 Suppose $r'_e=k$ so that $r'_i\le k$ for $i>e$. Then 
\begin{eqnarray*}
t_{ef}(R'',S'')&=& t_{ef}(R,S)+f-\sum_{i=e+1}^mr'_i\\
&\ge & \overline{\rho_k}  -ke-f+f-\sum_{i=e+1}^mr'_i\\
&=&\overline{\rho_k}-ke-\sum_{i=e+1}^m r'_i\\
&=&\overline{\rho_k}-\overline{\rho_k}\\
&=& 0.\\
\end{eqnarray*}
\end{proof}

The next theorem is the main result in this section.
In the statement of the theorem we use the integers $p_t,p_{t-1},\ldots,p_1$ as defined in (\ref{eq:last}).

\begin{theorem}\label{th:joint}
Let $R$ and $S$ be monotone nonincreasing, nonnegative integral vectors such that ${\mathcal A}(R,S)$ is nonempty. Let $t$ be a positive integer such that $\overline{\rho_t}=\overline{\rho_t}(R,S)\le n$. 
Then there exists a matrix $A\in {\mathcal A}(R,S)$ and a matrix $B\le A$ such that
\begin{enumerate}
\item[\rm (i)] $\rho_k(A)=\overline{\rho_k}$ for $k=1,2,\ldots,t$.
\item[\rm (ii)] $B$ contains exactly $\overline{\rho_t}$ $1$s,  where 
\begin{itemize}
\item[\rm (a)]
these $1$s  lie within the leading $\overline{\rho_1}\times \overline{\rho_t}$ submatrix of $B$,
\item[\rm (b)]  each of the first $\overline{\rho_t}$ columns of $B$ contains exactly one $1$, 
\item[\rm (c)] each of the first $\overline{\rho_1}$ rows of $B$ contains at most $t$ $1$s, with the first $p_t$ rows each containing $t$ $1$s, the next $p_{t-1}$ rows each containing $t-1$ $1$s, $\ldots$, the next $p_1$ rows each containing one $1$. 
\end{itemize}
\end{enumerate}
\end{theorem}

\begin{proof}
The theorem is an immediate  consequence of Lemma \ref{lem:joint}.
Note that $p_t$ is the maximum number of rows with $t$ $1$s we could have in such a $B$; then $p_{t-1}$ is the maximum number of
remaining rows we could have with $t-1$ $1$s, etc.
\end{proof}


\begin{thebibliography}{99}
 
\bibitem{RAB} R.A.~Brualdi, Transversal Matroids,  Chapter 5 in {\it Combinatorial Geometries}, N.~White ed., Cambridge Univ. Press, Cambridge, 1987, 72--97.



\bibitem{B} R.A.~Brualdi, {\it Combinatorial Matrix Classes}, Cambridge University Press, Cambridge, 2008.

\bibitem{BKMS} R.A.~Brualdi, K.P.~Kiernan, S.A.~Meyer, and M.W.~Schroeder, Combinatorial batch codes and transversal matroids,
{\it Adv. Math. Communications}, 4 (2010), 419-431.

\bibitem{BR} R.A.~Brualdi and H.J. Ryser, {\it Combinatorial Matrix Theory}, Cambridge University Press, Cambridge, 1991.

\bibitem{IKOS} Y.~Ishai, E.~Kushilevitz, R.~Ostrovsky, A.~Sahai,
Batch codes and their applications, in {\it Proceedings of the 36th Annual ACM Symposium on Theory of Computing}, ACM Press, New York, 262--271.

\bibitem{Ox} J.G.~Oxley, {\it Matroid Theory}, The Clarendon Press, Oxford University Press, New York, 1992.

\bibitem{PSW} M.B.~Paterson, D.R. Stinson, R.~Wei, Combinatorial batch codes, {\it Adv. Math. Communications}, 3 (2009), 13--27.

\end{thebibliography}
\end{document}